\documentclass[a4paper,11pt]{amsart}
\usepackage{amssymb}
\usepackage{ifthen}
\usepackage{graphicx}

\newtheorem{thm}{Theorem}
\newtheorem{cor}{Corollary}
\newtheorem{lem}{Lemma}

\newtheorem{rem}{Remark}
\newtheorem{prob}{Problem}
\newtheorem{conj}{Conjecture}
\theoremstyle{definition}
\newtheorem{example}[equation]{Example}

\newcounter {own}
\def\theown {\thesection       .\arabic{own}}

\newenvironment{pf}[1][]{%
 \vskip 3mm
 \noindent
 \ifthenelse{\equal{#1}{}}%
  {{\slshape Proof. }}%
  {{\slshape #1.} }%
 }%
{\qed\bigskip}

\newcounter{alphabet}
\newcounter{tmp}

\makeatletter
\newcommand{\Ref}[1]{\@ifundefined{r@#1}{}{\setcounter{tmp}{\ref{#1}}\Alph{tmp}}}
\makeatother

\newcounter{minutes}
\divide\time by 60
\newcounter{hours}
\multiply\time by 60
\addtocounter{minutes}{-\time}

\newcommand{\IC}{{\mathbb C}}
\newcommand{\ID}{{\mathbb D}}

\newcommand{\IZ}{{\mathbb Z}}

\def\be{\begin{equation}}
\def\ee{\end{equation}}

\newcommand{\bee}{\begin{enumerate}}
\newcommand{\eee}{\end{enumerate}}

\newcommand{\blem}{\begin{lem}}
\newcommand{\elem}{\end{lem}}
\newcommand{\bthm}{\begin{thm}}
\newcommand{\ethm}{\end{thm}}
\newcommand{\bcor}{\begin{cor}}
\newcommand{\ecor}{\end{cor}}
\newcommand{\beg}{\begin{example}}
\newcommand{\eeg}{\end{example}}
\newcommand{\begs}{\begin{examples}}
\newcommand{\eegs}{\end{examples}}
\newcommand{\bdefe}{\begin{defin}}
\newcommand{\edefe}{\end{defin}}
\newcommand{\bprob}{\begin{prob}}
\newcommand{\eprob}{\end{prob}}
\newcommand{\bei}{\begin{itemize}}
\newcommand{\eei}{\end{itemize}}

\newcommand{\bcon}{\begin{conj}}
\newcommand{\econ}{\end{conj}}
\newcommand{\bcons}{\begin{conjs}}
\newcommand{\econs}{\end{conjs}}
\newcommand{\bprop}{\begin{propo}}
\newcommand{\eprop}{\end{propo}}
\newcommand{\br}{\begin{rem}}
\newcommand{\er}{\end{rem}}
\newcommand{\brs}{\begin{rems}}
\newcommand{\ers}{\end{rems}}
\newcommand{\bo}{\begin{obser}}
\newcommand{\eo}{\end{obser}}
\newcommand{\bos}{\begin{obsers}}
\newcommand{\eos}{\end{obsers}}
\newcommand{\bpf}{\begin{pf}}
\newcommand{\epf}{\end{pf}}
\newcommand{\ba}{\begin{array}}
\newcommand{\ea}{\end{array}}
\newcommand{\beq}{\begin{eqnarray}}
\newcommand{\beqq}{\begin{eqnarray*}}
\newcommand{\eeq}{\end{eqnarray}}
\newcommand{\eeqq}{\end{eqnarray*}}

\newcommand{\ds}{\displaystyle}

\def\cc{\setcounter{equation}{0}   
\setcounter{figure}{0}\setcounter{table}{0}}

\begin{document}

\bibliographystyle{amsplain}

%

\title{The radius of univalence of the reciprocal of a product of two analytic functions}

\thanks{
File:~\jobname .tex,
          printed: \number\day-\number\month-\number\year,
          \thehours.\ifnum\theminutes<10{0}\fi\theminutes}
\author[\'A. Baricz]{\'Arp\'ad Baricz}
\address{\'A. Baricz, Department of Economics,  Babe\c{s}-Bolyai University, Cluj-Napoca 400591, Romania}
\address{Institute of Applied Mathematics, \'Obuda University, 1034 Budapest, Hungary}
\email{bariczocsi@yahoo.com}

\author[M. Obradovi\'{c}]{Milutin Obradovi\'{c}}
\address{M. Obradovi\'{c},
Department of Mathematics,
Faculty of Civil Engineering, University of Belgrade,
Bulevar Kralja Aleksandra 73, 11000
Belgrade, Serbia. }
\email{obrad@grf.bg.ac.rs}

\author[S. Ponnusamy]{Saminathan Ponnusamy
}
\address{S. Ponnusamy, Indian Statistical Institute (ISI), Chennai Centre, SETS (Society
for Electronic Transactions and security), MGR Knowledge City, CIT
Campus, Taramani, Chennai 600 113, India.}
\email{samy@isichennai.res.in, samy@iitm.ac.in}


\begin{abstract}
Let ${\mathcal A}$ denote the family of all   functions $f$ analytic
in the open unit disk $\ID$ with the normalization $f(0)=0= f'(0)-1$ and
${\mathcal S}$  be the class of univalent functions from ${\mathcal A}$.
In this paper,  we consider radius of univalence of $F$ defined by
$F(z)=z^{3}/(f(z)g(z))$, where $f$ and $g$ belong to some subclasses of ${\mathcal A}$ (for which
$f(z)/z$ and $g(z)/z$ are non-vanishing in $\ID$) and, in some cases in precise form, belonging to
some subclasses of ${\mathcal S}$. All the results are proved to be sharp. Applications of our
investigation through Bessel functions are also presented.
\end{abstract}


\maketitle
\pagestyle{myheadings}
\markboth{\'A. Baricz, M. Obradovi\'{c} and  S. Ponnusamy }{The radius of univalence of the reciprocal of a product of two analytic functions}
\cc

\section{Introduction and Main Results}
Let  $\ID$ be the open unit disk $\{z:\, |z|<1\}$ in the complex plane $\IC$, and $\mathcal H$ denote the linear
space of  analytic functions on $\ID$, with the topology of uniform convergence on compact subsets of $\ID$ so that
$\mathcal H$ is metrizable. In  $\mathcal H$, we consider the sub-collection
${\mathcal A}$  of functions $f\in {\mathcal H}$ with the normalization $f(0)=0= f'(0)-1$, and let
${\mathcal S}=\{f\in {\mathcal A}: \, \mbox{$f$ is univalent in $\ID$}\}.$ A number of geometric subclasses of
${\mathcal S}$ still enjoys the attention of many mathematicians in solving extremal
problems. We begin by recalling the standard ones:
By ${\mathcal C}$ and ${\mathcal S}^{\star}$
we denote the subclasses of ${\mathcal S}$ which consist of  convex (i.e. $f(\ID)$ is a convex domain) and
starlike functions (i.e. $f(\ID)$ is a domain starlike with respect to the origin), respectively.
Given a convex function $g\in\mathcal{C}$, we say that $f\in {\mathcal A}$ is said to belong to the class
$f\in {\mathcal K}_g$ if it satisfies the condition
\be\label{eq-1}
 {\rm Re\,}\left ( \frac{f'(z)}{g'(z)} \right )>0,
\quad \mbox{$z\in \ID$}.
\ee
Furthermore, a function $f\in  \mathcal{A}$ is said to
belong to ${\mathcal K}$ if $f\in {\mathcal K}_g$ for some convex function $g$ and thus,
$\mathcal{K}=\mathop{\cup}_{g\in\mathcal{C}}\mathcal{K}_g.$
Functions in $\mathcal{K}$ are known to be close-to-convex in $\ID$ (i.e. complement of $f(\ID)$ is connected by non-intersecting half lines),
and hence they belong to the class $\mathcal S$. In general, $g$ is not necessarily normalized.
Some particular choices of $g\in{\mathcal C}$ in ${\mathcal K}_g$ have
special role in many different contexts, for example, in deriving sufficient conditions for
functions to be in ${\mathcal K}$ (see \cite{MacG69}).

Now, we may recall the family ${\mathcal U}$ which has been studied in the recent years together
with its extension (see for instance, \cite{Aks58,FR-2006,OP-01,obpo-2007a,OPSV}):
$${\mathcal U} =  \{ f \in {\mathcal A} :\, \left |U_{f}(z)\right | < 1
~\mbox{ for $z\in \ID$}   \},
$$
where
\be\label{3-10eq1}
U_{f}(z)=f'(z)\left (\frac{z}{f(z)} \right )^{2}-1 .
\ee
The boundedness of $f'(z)(z/f(z))^2$ forces $f\in {\mathcal U}$ to be non-vanishing in
the punctured unit disk $0<|z|<1$. Hence $f'(z)\neq 0$ in $\ID$ and thus, $f$ is locally univalent in $\ID$.
Moreover, it is known \cite{Aks58} that $f\in  {\mathcal U}$ univalent in $\ID$ and thus, ${\mathcal U}\subset {\mathcal S}$.
It is known \cite{FR-2006,OP-01,PV2005} that neither ${\mathcal U}$ is included in ${\mathcal S}^{\star}$ nor
includes ${\mathcal S}^{\star}$ (see \cite{FR-2006}). In fact, $\mathcal{U}$ is not a subset of
$\mathcal{S}^{*}$ as the function
$$f_1(z)=\frac{z}{1+\frac{1}{2}z+\frac{1}{2}z^3}
$$
demonstrates.  In view of this reasoning, general classes of functions was
considered by Obradovi\'{c} and Ponnusamy, for example,   for $0<\lambda \leq 1$, the class
$${\mathcal U}(\lambda ) =  \{ f \in {\mathcal A} :\, \left |U_{f}(z)\right | < \lambda
~\mbox{ for $z\in \ID$} \}.
$$
In order to compare with the area principle of Gronwall for meromorphic functions (see \cite[p.~18, Theorem 1.3]{Pomm}
and \cite[p.~29, Theorem 2.1]{Du}), we observe that mappings $f\in \mathcal{S}$
can be associated with the mappings $F\in \Sigma$, namely univalent functions $F$ of the form,
$$ F(\zeta )=\zeta +\sum_{n=0}^{\infty}c_n \zeta ^{-n}, \quad |\zeta | >1,
$$
which satisfies the condition $F(\zeta ) \neq 0$ for $|\zeta|>1$, by the correspondence
$$F(\zeta ) =\frac{1}{f(1/\zeta)}, \quad |\zeta |>1 .
$$
Using the change of variable $\zeta =1/z$, the association $f(z)=1/F(1/z)$ quickly yields the formula
$$F'(\zeta)-1=U_f(z),
$$
where $U_f$ is defined by \eqref{3-10eq1} and $a_2=f''(0)/2=-c_0$. The last relation provides the close connection that exists between
$f\in \mathcal U$ and  $F\in \Sigma$ of meromorphic univalent functions $F(\zeta)$ satisfying the condition $|F'(\zeta)-1|<1$ for $ |\zeta | >1$.
One of the interesting observations is that each function in
$${\mathcal S}_{\IZ}= \left \{z, ~~\frac{z}{(1\pm z)^2},~~\frac{z}{1\pm z},~~\frac{z}{1\pm z^2}, ~~ \frac{z}{1\pm z+z^2}\right \}
$$
belongs to $\mathcal{U}$. Also, it is well-known that functions in ${\mathcal S}_{\IZ}$
are the only functions in $\mathcal S$ having integral coefficients in the Taylor series
expansions of $f\in {\mathcal S}$.   Moreover, every $f\in \mathcal{U}(\lambda )$ can be expressed as (cf. \cite{obpo-2007a})
$$\frac{z}{f(z)}=1-a_{2}z-\lambda z\int_{0}^{z}\frac{\omega(t)}{t^{2}}dt, ~ a_2 =\frac{f''(0)}{2},
$$
for some  $\omega\in{\mathcal B}_1$, where ${\mathcal B}_1$ denotes the
class of analytic functions in the unit disk $\ID$ such  $\omega(0)=\omega '(0)=0$  and $|\omega(z)|<1$ for $z\in\ID$.
More recently, Vasudev and Yanagihara  \cite{VY2013} discussed the class ${\mathcal U}(\lambda )$ in geometric perspectives.

When we say that  $f\in\mathcal{U}$ in $|z|<r$ it means that the condition  $|U_{f}(z)| < 1$ holds in the sub disk $|z|<r$ instead of
the full unit disk $\ID$, which is indeed same as saying that $r^{-1}f(rz)$ belongs to the class $\mathcal{U}$. It is convenient to use the
same convention at similar situations. Some related classes to these families may now be recalled:
\beqq
\mathcal{R} &=& \{f\in{\mathcal A}:\, {\rm Re\,}f'(z)>0~\mbox{ for $z\in \ID$}\}\\
{\mathcal P}(1/2) &=& \left \{f\in {\mathcal A}:\, {\rm Re\,}(f(z)/z)>1/2 ~\mbox{ for $z\in \ID$}\right \} ~\mbox{ and }\\
{\mathcal C}(-1/2) &=& \left \{f\in {\mathcal A}:\,{\rm Re\,}\left(1+\frac{zf''(z)}{f'(z)}\right)>-\frac{1}{2}
~\mbox{ for $z\in \ID$} \right \}.
\eeqq
We remark that ${\mathcal C}\subset{\mathcal S}^{\star}(1/2)  \subset {\mathcal P}(1/2)$. Here
${\mathcal S}^{\star}(1/2)$ is the class of starlike functions of order $1/2$.
Analytically,
$f\in{\mathcal S}^{\star}(1/2)$ if and only if $f\in  \mathcal{A}$ and ${\rm Re\,} (zf'(z)/f(z))>1/2$ in $\ID$.
Also, it is well-known that (see \cite{Um52}) each function in ${\mathcal C}(-1/2)$ is indeed convex in some direction
and hence, ${\mathcal C}(-1/2)$ is included in ${\mathcal K}$.
Many properties of these classes and their generalizations
have been studied extensively in the literature (see \cite{Du,Go,Pomm}). There are many necessary and sufficient coefficient
conditions for functions to be in these classes.  The following well-known necessary conditions for functions in
$\mathcal{S}$ and some of its subclasses are needed for our investigation:
\bee
\item[(a)] if $f\in \mathcal{R}$, then $|a_n|\leq 2/n$ for all $n\geq 2$
\item[(b)] if $f\in {\mathcal P}(1/2)\cup {\mathcal C}\cup {\mathcal S}^{\star}(1/2)$, then $|a_n|\leq 1$ for all $n\geq 2$
\item[(c)] if $f\in \mathcal{S}$, then $|a_n|\leq n$ for all $n\geq 2$
\item[(d)] if $f\in {\mathcal C}(-1/2)$, then $|a_n|\leq (n+1)/2$ for all $n\geq 2$
(see \cite[Corollary 2]{Shah73} and \cite{samy-hiroshi-swadesh}).
\eee

In this paper we are mainly interested in the following question:

\bprob\label{prob1}
Suppose that $f,g\in {\mathcal A}$ such that $f(z)/z$ and $g(z)/z$ are non-vanishing in $\ID$. If $F$ is of the form
\be\label{3-10eq2}
F(z)=\frac{z^{3}}{f(z)g(z)} ~\mbox{ for $z\in \ID$},
\ee
what can be said about the univalence of $F$?
\eprob

Clearly, $F$ could be considered as an operation acting on the space of analytic functions into another with some standard procedure.
Recall that if $f\in \mathcal{S}$ and is of the form  
$$\frac{z}{f(z)} = 1+\sum_{n=1}^{\infty}c_nz^n 
$$
then the well-known Area Theorem \cite[Theorem 11 on p.193 of Vol. 2]{Go} gives the following necessary coefficient condition
$$\sum_{n=2}^\infty (n-1)|c_n|^2\leq 1.
$$
This condition has been very helpful in solving many different problems in the theory of univalent functions.
However, a natural question concerning Problem \ref{prob1} is the following: How do the Taylor coefficients of $f$ and $g$ influence the
property of the normalized analytic function $F$? Do the restrictions on the Taylor coefficients of $f$ and $g$ provide examples of
functions from $\mathcal U$, for instance?  The answer is yes and we consider this problem for a few cases but from the
proof one can easily see that many general results could be obtained using our approach.

Throughout the discussion $f,g\in {\mathcal A}$ will be of the form
\be\label{3-10eq3}
f(z)=z+\sum_{n=2}^{\infty}a_{n}z^{n}~\mbox{ and }~g(z)=z+\sum_{n=2}^{\infty}b_{n}z^{n}.
\ee
The paper is organized as follows. In Section \ref{sec2}, we state our main results and their direct consequences. Section \ref{sec3}
outlines the basic idea of the proof and present a partial list of situations for which our conclusions in the main theorems continue
to hold with varying hypotheses. We present proofs of our main theorems in Section \ref{sec4}. In Section \ref{sec5}, we present a case where
our results can be applied for Bessel functions of the first kind. The concluding section presents further remarks on our
methodology.

\section{Main Results and their consequences}\label{sec2}
We begin to present a first set of results based on the restrictions on the coefficients of $f$ and $g$.

\bthm\label{3-13th1}
Let $f,g,F\in {\mathcal A}$ be as in Problem \ref{prob1}. Then we have the following:
\begin{enumerate}
\item[{\rm \textbf{(a)}}]  If $|a_{n}|\leq 1$ and $|b_{n}|\leq 1$ for $n\geq 2$,
then $F\in \mathcal{U}$ in the disk $|z|<r_1=1/3$.

\item[{\rm \textbf{(b)}}]  If $|a_{n}|\leq n$ and $|b_{n}|\leq 1$ for $n\geq 2$,
then $F\in \mathcal{U}$ in the disk $|z|<r_2=1/4.$

\item[{\rm \textbf{(c)}}]   If $|a_{n}|\leq n$ and $|b_{n}|\leq n$ for $n\geq 2$,
then $F\in \mathcal{U}$ in the disk $|z|<r_3=1/5.$
\end{enumerate}
All these results are best possible (as for the univalence is concerned).
\ethm

The following corollary is a consequence of Theorem \ref{3-13th1} and so, we do not need a detailed explanation but
it is sufficient to recall that if $f\in {\mathcal P}(1/2)\cup {\mathcal C}\cup {\mathcal S}^{\star}(1/2)$ and
$g\in {\mathcal S}$, then we have $|a_n|\leq 1$ and $|b_n|\leq n$ for $n\geq 2$, respectively.

\bcor
Let $F$ be defined by \eqref{3-10eq2}. Then we have the following:
\begin{enumerate}
\item[{\rm \textbf{(a)}}] $f,g\in {\mathcal P}(1/2)\cup {\mathcal C}\cup {\mathcal S}^{\star}(1/2)$
imply that $F\in \mathcal{U}$ in the disk $|z|<1/3.$ As for the univalence, the radius $1/3$ cannot be replaced by
a bigger number and thus the result is sharp.

\item[{\rm \textbf{(b)}}] If $f\in {\mathcal S}$ and $g\in {\mathcal P}(1/2)\cup {\mathcal C} \cup {\mathcal S}^{\star}(1/2)$,
then $F\in \mathcal{U}$ in the disk $|z|<1/4$, and the radius $1/4$ is best possible (as for the univalence).

\item[{\rm \textbf{(c)}}]  If $f,g\in {\mathcal S}$, then $F\in \mathcal{U}$ in the disk $|z|<1/5$
and the radius $1/5$ is sharp.
\end{enumerate}
\ecor

We now state our next result.

\bthm\label{3-13th1a}
Let $f,g,F\in {\mathcal A}$ be as in Problem \ref{prob1}. Then we have the following:
\begin{enumerate}
\item[{\rm \textbf{(a)}}]  If $|a_{n}|\leq 1$ and $|b_{n}|\leq 2/n$ for $n\geq 2$,
then $F\in \mathcal{U}$ in the disk $|z|<r_4 $, where $r_4\approx 0.36027 $ is the root of the equation
\be\label{eq10f}
2r-3 +\frac{2(3r-2)}{r}\log (1-r)=0
\ee
in the interval $(0,1)$.

\item[{\rm \textbf{(b)}}]  If $|a_{n}|\leq n$ and $|b_{n}|\leq 2/n$ for $n\geq 2$,
then $F\in \mathcal{U}$ in the disk $|z|<r_5 $, where $r_5\approx 0.26073$ is the root of the equation
\be\label{eq10y}
3(1-r) -\frac{4(2r-1)}{r}\log (1-r)=0
\ee
in the interval $(0,1)$.

\item[{\rm \textbf{(c)}}]   If $|a_{n}|\leq 2/n$ and $|b_{n}|\leq 2/n$ for $n\geq 2$,
then $F\in \mathcal{U}$ in the disk $|z|<r_6 $, where $r_6\approx 0.399185$ is the unique root of the equation
\be\label{eq10d}
\frac{5-r}{1-r} +\frac{6}{r}\log (1-r)=0
\ee
in the interval $(0,1)$.
\end{enumerate}
All these results are best possible (as for the univalence is concerned).
\ethm


\bcor
Let $F$ be defined by \eqref{3-10eq2} and $r_j$ $(j=4,5,6)$ be as in Theorem \ref{3-13th1a}. Then we have the following:
\begin{enumerate}
\item[{\rm \textbf{(a)}}] $f\in {\mathcal P}(1/2)\cup {\mathcal C}\cup {\mathcal S}^{\star}(1/2)$  and  $g\in \mathcal{R}$
imply that $F\in \mathcal{U}$ in the disk $|z|<r_4.$ As for the univalence the result is sharp.

\item[{\rm \textbf{(b)}}] If $f\in {\mathcal S}$ and $g\in \mathcal{R}$,
then $F\in \mathcal{U}$ in the disk $|z|<r_5$, and the result is best possible (as for the univalence).

\item[{\rm \textbf{(c)}}]  If $f,g\in \mathcal{R}$, then $F\in \mathcal{U}$ in the disk $|z|<r_6.$
The result is sharp.
\end{enumerate}
\ecor

Our final result follows.

\bthm\label{3-13th2}
Let $f,g,F\in {\mathcal A}$ be as in Problem \ref{prob1}. Then we have the following:
\begin{enumerate}
\item[{\rm \textbf{(a)}}]  If $|a_{n}|\leq (n+1)/2$ and $|b_{n}|\leq 1$ for $n\geq 2$,
then $F\in \mathcal{U}$ in the disk $|z|<r_7=\frac{4-\sqrt{10}}{3} \approx 0.27924$.
As for the univalence the result is sharp.

\item[{\rm \textbf{(b)}}] If  $|a_{n}|\leq (n+1)/2$ and $|b_{n}|\leq n$ for $n\geq 2$,
then $F\in \mathcal{U}$ in the disk $|z|<r_8=\frac{5-\sqrt{17}}{4}\approx 0.21922$
and the result is best possible (as for the univalence).

\item[{\rm \textbf{(c)}}] If  $|a_{n}|\leq (n+1)/2$ and $|b_{n}|\leq 2/n$ for $n\geq 2$,
then $F\in \mathcal{U}$ in the disk $|z|<r_9$, where $r_9\approx  0.29399$
is the unique root of the equation
\be\label{eq10e}
(r-3)(1-r) -\frac{4-9r+3r^2}{r}\log (1-r)=0
\ee
in the interval $(0,1)$. The result is best possible (as for the univalence).
\end{enumerate}
\ethm

We end the section with a corollary which is a consequence of Theorem \ref{3-13th2}.

\bcor
Let the function $F$ be defined with \eqref{3-10eq2}  and $r_j$ $(j=7,8,9)$ be as in Theorem \ref{3-13th2}.  Then we have the following:
\begin{enumerate}
\item[{\rm \textbf{(a)}}] If $f\in{\mathcal C}(-1/2)$ and  $g\in  {\mathcal P}(1/2)\cup {\mathcal C}\cup {\mathcal S}^{\star}(1/2)$,
then $F\in \mathcal{U}$ in the disk $|z|<r_7 .$ As for the univalence the result is sharp.

\item[{\rm \textbf{(b)}}] If $f\in{\mathcal C}(-1/2)$ and  $g\in {\mathcal S}$,
then $F\in \mathcal{U}$ in the disk $|z|<r_8$
and the result is best possible (as for the univalence).

\item[{\rm \textbf{(c)}}]  If $f\in{\mathcal C}(-1/2)$ and $g\in \mathcal{R}$, then $F\in \mathcal{U}$ in the disk $|z|<r_9.$
The result is sharp.
\end{enumerate}
\ecor

\section{General Guidelines for the Proofs of our results}\label{sec3}
\subsection{Basic idea of the Proof}
Assume the hypothesis of Problem \ref{prob1} on $f,g,F\in {\mathcal A}$  and the power series representations for $f$ and $g$
given by \eqref{3-10eq3}. Using the notations for $f$ and $F$ given by \eqref{3-10eq1} and \eqref{3-10eq2}, respectively,  a computation gives
\beqq
U_{F}(z)&=& -z^{2}\left(\frac{1}{F(z)}-\frac{1}{z}\right)'\\
&=& -z^{2}\left(\frac{1}{z}\left[\frac{f(z)}{z}\frac{g(z)}{z}-1 \right]\right)' \\
&=&\frac{f(z)}{z}\frac{g(z)}{z}-1-z\left(\frac{f(z)}{z}\right)'\frac{g(z)}{z}-
z\left(\frac{g(z)}{z}\right)'\frac{f(z)}{z}
\eeqq
which may be further decomposed as
\beqq
U_{F}(z) &=& \frac{g(z)}{z}\left[\frac{f(z)}{z}-z\left(\frac{f(z)}{z}\right)'-1 \right]
+\frac{f(z)}{z}\left[\frac{g(z)}{z}-z\left(\frac{g(z)}{z}\right)'-1 \right]\\
&& -\left(\frac{f(z)}{z}-1\right)\left(\frac{g(z)}{z}-1\right)
\eeqq
and by the triangle inequality, we have
\beq\label{3-10eq4}
\left|U_{F}(z)\right|&\leq&  \left| \frac{g(z)}{z}\right|\left|\frac{f(z)}{z}-z\left(\frac{f(z)}{z}\right)'-1 \right|
+\left|\frac{f(z)}{z}\right|\left|\frac{g(z)}{z}-z\left(\frac{g(z)}{z}\right)'-1 \right|\nonumber\\
&& + \left|\frac{f(z)}{z}-1\right|\left|\frac{g(z)}{z}-1\right|.
\eeq
The aim is to find estimate for $|U_{F}(z)|$ when some properties about $f$ and $g$ or some appropriate coefficients conditions
for $a_n$ and $b_n$ are given. When the coefficient conditions are given as in Theorems \ref{3-13th1}, \ref{3-13th1a} and \ref{3-13th2},
we may proceed as follows: We use the representation of $f$ in \eqref{3-10eq3} and obtain for $|z|\leq r$,
\be\label{eq2}
\left|\frac{f(z)}{z}-1\right|\leq \sum_{n=2}^{\infty}|a_{n}|r^{n-1} ~\mbox{ and }~
\left|\frac{f(z)}{z}\right| \leq 1+\sum_{n=2}^{\infty}|a_{n}|r^{n-1}
\ee
and similarly,
\be\label{eq3}
\left|\frac{f(z)}{z}-z\left(\frac{f(z)}{z}\right)'-1 \right|=\left|-\sum_{n=3}^{\infty}(n-2)a_{n}z^{n-1}\right|
\leq \sum_{n=3}^{\infty}(n-2)|a_{n}|r^{n-1}.
\ee
In general, for many basic subclasses of functions $f\in{\mathcal S}$, a necessary coefficient condition will be of the form
\be\label{eq4}
|a_n|\leq A_n(a,b,c)=a +(b/n)+cn ~\mbox{ for all $n\geq 2$},
\ee
where $a,b,c$ are suitable real numbers such that $A_n(a,b,c)>0$ for  all $n\geq 2$.
This is our basic here and in Section \ref{sec6}, we include a partial list of these cases.

\subsection{Basic inequalities for the proof of our main results}
The following equalities will be used in the sequel:
\be\label{eq5}
\frac{-\log (1-r)}{r} =\sum_{n=1}^{\infty}\frac{r^{n-1}}{n},~~\frac{1}{(1-r)^2} =\sum_{n=1}^{\infty}nr^{n-1},~~
\frac{1+r}{(1-r)^3} =\sum_{n=1}^{\infty}n^2r^{n-1}.
\ee
If $a_n$ satisfies the condition \eqref{eq4}, then with the help of \eqref{eq5} the first relation in \eqref{eq2} gives
\beqq
\left|\frac{f(z)}{z}-1\right|&\leq& a\sum_{n=2}^{\infty}r^{n-1} +b\sum_{n=2}^{\infty}\frac{r^{n-1}}{n} +c\sum_{n=2}^{\infty}nr^{n-1}\\
&=& a\left [ \frac{r}{1-r}\right ]+ b\left [ -\frac{\log (1-r)}{r} -1\right ]+c \left [\frac{1}{(1-r)^2}-1 \right ]
\eeqq
which implies that
\be\label{eq6}
\left|\frac{f(z)}{z}-1\right|\leq -(a+b+c)+ \frac{a}{1-r} -\frac{b\log (1-r)}{r} +\frac{c}{(1-r)^2}
\ee
and thus,
\be\label{eq7}
\left|\frac{f(z)}{z}\right|\leq 1-(a+b+c)+ \frac{a}{1-r} -\frac{b\log (1-r)}{r} +\frac{c}{(1-r)^2}.
\ee
Similarly, using the condition \eqref{eq4} and the rearrangement
$$(n-2)[a +(b/n)+cn]= -2a+b+(a-2c)n+cn^2-\frac{2b}{n},
$$
the relation \eqref{eq3} becomes

\vspace{8pt}

$\ds \left|\frac{f(z)}{z}-z\left(\frac{f(z)}{z}\right)'-1 \right|$
\beqq
 &\leq& \sum_{n=2}^{\infty}(n-2) [a +(b/n)+cn]r^{n-1}\\
&=& (-2a+b)\sum_{n=2}^{\infty}r^{n-1}+ (a-2c)\sum_{n=2}^{\infty}nr^{n-1} +c\sum_{n=2}^{\infty}n^2r^{n-1} -2b\sum_{n=2}^{\infty}\frac{r^{n-1}}{n}\\
&=& (-2a+b)\left [ \frac{1}{1-r}-1\right ] +  (a-2c)\left [\frac{1}{(1-r)^2}-1 \right ] +c \left [\frac{1+r}{(1-r)^3}-1 \right ] \\
&& ~~~+ 2b\left [ \frac{\log (1-r)}{r} +1\right ]\\
\eeqq
and thus, a simplification gives
\be\label{eq8}
\left|\frac{f(z)}{z}-z\left(\frac{f(z)}{z}\right)'-1 \right|\leq
a+b+c- \frac{2a-b}{1-r}   +\frac{a-2c}{(1-r)^2} +\frac{c(1+r)}{(1-r)^3}+\frac{2b\log (1-r)}{r}.
\ee

\section{Proofs of main theorems}\label{sec4}

From the proof of our theorems, one can obtain a number of general theorems for various other choices for $a,b,c$.
A similar approach is possible to various other situations although we omit them here.

\subsection{Proof of Theorem \ref{3-13th1}}
\textbf{(a)} Suppose that $|a_{n}|\leq 1$ and $|b_{n}|\leq 1$ for $n\geq 2$. This means that we choose $a=1$ and $b=c=0$ in
\eqref{eq4} and thus,  by using the relations \eqref{eq6}, \eqref{eq7} and \eqref{eq8} (with $a=1$, $b=c=0$), one has
\be\label{eq9}
\left \{   \ba{l} \ds \left|\frac{f(z)}{z}-1\right|\leq \frac{r}{1-r}, ~~ \left|\frac{f(z)}{z}\right|\leq \frac{1}{1-r}, ~\mbox{ and }\\[4mm]
\ds \left|\frac{f(z)}{z}-z\left(\frac{f(z)}{z}\right)'-1 \right|\leq \frac{r^2}{(1-r)^2}
\ea\right .
\ee
for $|z|\leq r$.
Similar inequalities hold for $g$ (because $|b_{n}|\leq 1$ for $n\geq 2$) and thus, using these bounds,  \eqref{3-10eq4} gives us the estimate
$$\left|U_{F}(z)\right|\leq \frac{2r^2}{(1-r)^{3}}+ \frac{r^2}{(1-r)^{2}} =\frac{3r-1}{(1-r)^{3}}+1
$$
which is clearly less than $1$ if  $0\leq r<r_1=1/3$, i.e. $F\in {\mathcal U}$ for $|z|<1/3$. In particular, $F$ is univalent
for $|z|<{1}/{3}$. If we choose $f(z)=g(z)=z/(1-z)$, then we obtain $F(z)=z(1-z)^{2}$ and $F'(z)=1-4z+3z^{2}$
showing that  $F'(1/3)=0$ and thus,  $F$ is not univalent in any bigger disk. Thus, the radius $1/3$ is best possible.

\textbf{(b)} Suppose that $|a_{n}|\leq n$ and $|b_{n}|\leq 1$ for $n\geq 2$. As in Case \textbf{(a)}
but by using the relations \eqref{eq6}, \eqref{eq7} and \eqref{eq8} with $a=b=0$ and $c=1$, one has for $f$ the inequalities
\be\label{eq10}
\left \{   \ba{l} \ds \left|\frac{f(z)}{z}-1\right|\leq \frac{r(2-r)}{(1-r)^2}, ~~ \left|\frac{f(z)}{z}\right|\leq \frac{1}{(1-r)^2}, ~\mbox{ and }\\[4mm]
\ds \left|\frac{f(z)}{z}-z\left(\frac{f(z)}{z}\right)'-1 \right|\leq \frac{r^2(3-r)}{(1-r)^3}
\ea\right .
\ee
for $|z|\leq r$.
Inequalities of the form \eqref{eq9} holds for the function $g$ in place of $f$. This observation and \eqref{eq10} give the
following estimate for \eqref{3-10eq4}:
$$\left|U_{F}(z)\right|\leq \frac{r^2(3-r)}{(1-r)^{4}}+ \frac{r^2}{(1-r)^{4}} +\frac{r^2(2-r)}{(1-r)^{3}} =\frac{4r-1}{(1-r)^{4}}+1
$$
which implies that $\left|U_{F}(z)\right|<1$ for $0\leq r<r_2=1/4$, i.e. $F\in {\mathcal U}$ for $|z|<1/4$.
If we put $f(z)=z/(1-z)^{2}$ and $g(z)=z/(1-z)$, then we see that
$F(z)=z(1-z)^{3}$ and $F'(z)=(1-4z)(1-z)^{2}$. Thus,  $F'(1/4)=0$ and hence,
$F$ cannot be univalent in a larger disk. Consequently, the radius $1/4$ is best possible as for univalence of $F$ is concerned.

\textbf{(c)} Assume that $|a_{n}|\leq n$ and $|b_{n}|\leq n$ for $n\geq 2$. Then \eqref{eq10} holds and it also holds for $g$ in place of $f$.
In view of these observations, as in cases \textbf{(a)} and \textbf{(b)}, it follows easily that
$$\left|U_{F}(z)\right|\leq \frac{2r^2(3-r)}{(1-r)^{5}}+ \frac{r^2(2-r)^{2}}{(1-r)^{4}} =
\frac{5r-1}{(1-r)^{5}}+1,
$$
and thus, $\left|U_{F}(z)\right|<1$ for $0\leq r<r_3=1/5$, i.e. $F\in {\mathcal U}$ for $|z|<1/5$.
By putting $f(z)=g(z)=z/(1-z)^{2}$, we get $F(z)=z(1-z)^{4}$ and so, $F'(z)=(1-5z)(1-z)^{3}$.
We see that $F'(1/5)=0$ which means that $F$ cannot be univalent in the disk $|z|<r$ if $r>1/5$. Again,
the radius $1/5$ is best possible.
\hfill $\Box$



\subsection{Proof of Theorem \ref{3-13th1a}}
\textbf{(a)} Suppose that $|a_{n}|\leq 1$ and $|b_{n}|\leq 2/n$ for $n\geq 2$. Then \eqref{eq9} holds. Since
$|b_{n}|\leq 2/n$ for $n\geq 2$ (i.e. considering the choice $a=c=0$ and $b=2$ in
\eqref{eq4}),  by using the relations \eqref{eq6}, \eqref{eq7} and \eqref{eq8}, one has for $g$ the inequalities
\be\label{eq10c}
\left \{   \ba{l} \ds \left|\frac{g(z)}{z}-1\right|\leq -2-\frac{2}{r}\log (1-r), ~~ \left|\frac{g(z)}{z}\right|\leq -1-\frac{2}{r}\log (1-r),
~\mbox{ and }\\[4mm]
\ds \left|\frac{g(z)}{z}-z\left(\frac{g(z)}{z}\right)'-1 \right|\leq  2+\frac{2}{1-r} +\frac{4}{r}\log (1-r) = \frac{4-2r}{1-r} +\frac{4}{r}\log (1-r)
\ea\right .
\ee
for $|z|\leq r$. Using \eqref{eq9} and \eqref{eq10c}, we have the estimate
\beqq
\left|U_{F}(z)\right|& \leq & \left ( -1-\frac{2}{r}\log (1-r) \right )\frac{r^2}{(1-r)^{2}}+
\frac{1}{1-r}\left ( \frac{4-2r}{1-r} +\frac{4}{r}\log (1-r) \right ) \\
& & +\frac{r}{1-r}\left ( -2-\frac{2}{r}\log (1-r) \right )\\
&=& 1- \frac{1}{(1-r)^{2}}\left ( 2r-3+\frac{2(3r-2)}{r}\log (1-r)\right )
\eeqq
which is less than $1$ provided $\phi _1(r)>0$, where
$$\phi _1(r)= 2r-3+\frac{2(3r-2)}{r}\log (1-r) .
$$
This gives the conclusion $\left|U_{F}(z)\right|<1$ for $0\leq r<r_4$, i.e. $F\in {\mathcal U}$ for $|z|<r_4$, where $r_4$ is the root
of the equation \eqref{eq10f}, namely, $\phi _1 (r)=0$. If we put $f(z)=z/(1-z)$ and $g(z)=-z -2\log (1-z)$, then we see that $a_n=1$ and
$b_n=2/n$ for $n\geq 2$. We see that
$$F(z)=\frac{z(1-z)}{-1- (2/z)\log (1-z)}
$$
and thus,
$$F'(z)=\frac{2z-3+\frac{2(3z-2)}{z}\log (1-z)}{[-1- (2/z)\log (1-z)]^{2}}.
$$
It follows that  $F'(r_4)=0$ and hence, $F$ cannot be univalent in a larger disk.

\textbf{(b)} Suppose that $|a_{n}|\leq n$ and $|b_{n}|\leq 2/n$ for $n\geq 2$. Then \eqref{eq10} and \eqref{eq10c} hold.
As a consequence of these observations, we have the estimate
\beqq
\left|U_{F}(z)\right|& \leq & \left ( -1-\frac{2}{r}\log (1-r) \right )\frac{r^2(3-r)}{(1-r)^{3}}+
\frac{1}{(1-r)^{2}}\left ( \frac{4-2r}{1-r} +\frac{4}{r}\log (1-r) \right ) \\
&& +\frac{r(2-r)}{(1-r)^{2}}\left ( -2-\frac{2}{r}\log (1-r) \right )\\
&=& 1-\frac{1}{(1-r)^{3}}\left ( \frac{4(2r-1)}{r}\log (1-r)-3(1-r) \right )
\eeqq
which is less than $1$ provided $\phi _2(r)>0$, where
$$\phi _2(r)= \frac{4(2r-1)}{r}\log (1-r)-3(1-r).
$$
This gives the conclusion $\left|U_{F}(z)\right|<1$ for $0\leq r<r_5$, i.e. $F\in {\mathcal U}$ for $|z|<r_5$, where $r_5$ is the root
of the equation \eqref{eq10y}, namely, $\phi _2(r)=0$.

If we put $f(z)=z/(1-z)^{2}$ and $g(z)=-z -2\log (1-z)$, then we see that $a_n=n$ and $b_n=2/n$ for $n\geq 2$. We see that
$$F(z)=\frac{z(1-z)^2}{-1- (2/z)\log (1-z)}
$$
and thus,
$$F'(z)=-\frac{(1-z)\left [3(1-z) -\frac{4(2z-1)}{z}\log (1-z)\right ]}{[-1- (2/z)\log (1-z)]^{2}}.
$$
It follows that  $F'(r_5)=0$ and hence, $F$ cannot be univalent in a larger disk.

\textbf{(c)} As in cases \textbf{(a)} and \textbf{(c)}  of Theorem \ref{3-13th1}, we first observe that \eqref{eq10c} holds (by hypotheses)
both for $f$ and $g$.
This observation  and a computation give the following estimate for \eqref{3-10eq4}:
\beqq
\left|U_{F}(z)\right|& \leq & 2\left ( -1-\frac{2}{r}\log (1-r) \right ) \left ( \frac{4-2r}{1-r} +\frac{4}{r}\log (1-r) \right ) +\left ( -2-\frac{2}{r}\log (1-r) \right )^2 \\
&=& 1 +\left ( -1-\frac{2}{r}\log (1-r) \right ) \left [ \frac{5-r}{1-r} +\frac{6}{r}\log (1-r) \right ]
\eeqq
which implies that $\left|U_{F}(z)\right|<1$ for $0\leq r<r_6$, i.e. $F\in {\mathcal U}$ for $|z|<r_6$, where $r_6$ is the root
of the equation \eqref{eq10d}, namely, $\phi _3(r)=0$, where
$$\phi _3(r)= \frac{5-r}{1-r} +\frac{6}{r}\log (1-r).
$$

If we put $f(z)=g(z)=-z -2\log (1-z)$, then we see that $a_n=b_n=2/n$ for $n\geq 2$, $F(z)=z[-1- (2/z)\log (1-z)]^{-2}$ and thus
$$F'(z)=-\frac{\frac{5-z}{1-z} +\frac{6}{z}\log (1-z) }{[-1- (2/z)\log (1-z)]^{3}}.
$$
It follows that  $F'(r_6)=0$ and hence, $F$ cannot be univalent in a larger disk.
\hfill $\Box$


\subsection{Proof of Theorem \ref{3-13th2}}
\textbf{(a)} Let $f,g,F\in {\mathcal A}$ be as in Problem \ref{prob1},   $|a_{n}|\leq (n+1)/2$ and $|b_{n}|\leq 1$ for $n\geq 2$.
As before we see that \eqref{eq9} holds with $g$ in place of $f$. Moreover, since $|a_{n}|\leq (n+1)/2$, $f$ satisfies the inequalities
\be\label{eq10g}
\left \{   \ba{l} \ds \left|\frac{f(z)}{z}-1\right|\leq \frac{3r-2r^{2}}{2(1-r)^{2}}, ~~
\left|\frac{f(z)}{z}\right|\leq \frac{2-r}{2(1-r)^{2}}, ~\mbox{ and }\\[4mm]
\ds \left|\frac{f(z)}{z}-z\left(\frac{f(z)}{z}\right)'-1 \right|\leq \frac{2r^{2}-r^{3}}{(1-r)^{3}}.
\ea\right .
\ee
for $|z|\leq r$. As a consequence of these inequalities, the relation
\eqref{3-10eq4} simplifies to
\beqq
\left|U_{F}(z)\right|&\leq& \frac{1}{1-r}\left (\frac{2r^2-r^3}{(1-r)^{3}}\right )+
\frac{2-r}{2(1-r)^{2}}\left (\frac{r^2}{(1-r)^{2}}\right )+
\frac{r}{1-r}\left ( \frac{3r-2r^{2}}{2(1-r)^{2}}\right )\\
&=& \frac{9r^{2}-8r^{3}+2r^{4}}{(1-r)^{4}}\\
&=& 1-  \frac{3r^2-8r+2}{2(1-r)^{4}}
\eeqq
which is less than $1$ provided $3r^2-8r+2>0$, i.e.
if  $0\leq r<r_7=\frac{4-\sqrt{10}}{3}$. This means that $F\in {\mathcal U}$ in the disk $|z|<r_7$ and so, $F$ is univalent
in that disk.

As for the sharpness of the radius, we choose
$$f(z)=\frac{z(2-z)}{2(1-z)^{2}}=z+\sum_{n=2}^{\infty}\left (\frac{n+1}{2}\right )z^{n} ~ \mbox{ and }~g(z)=\frac{z}{1-z},
$$
so that $a_n= (n+1)/2$, $b_n=1$ and thus, $F$ defined by \eqref{3-10eq2} takes the form
$$ F(z)=\frac{2z(1-z)^{3}}{2-z}.
$$
We compute that
$$F'(z)=\frac{2(1-z)^{2}(3z^{2}-8z+2)}{\left(2-z\right)^{2}}
$$
and thus, $F'(r_{7})=0$ which means that $F$ is not univalent in a bigger disk.

\textbf{(b)} Let $f,g,F\in {\mathcal A}$ be as in Problem \ref{prob1}. For $f$ and $g$, by assumption, we have
$|a_{n}|\leq (n+1)/2$ and $|b_{n}|\leq n$ for $n\geq 2$, respectively.  Following the arguments as in the other cases,
we easily have the estimate
\beqq
\left|U_{F}(z)\right|&\leq& \frac{1}{(1-r)^{2}}\left (\frac{2r^2-r^3}{(1-r)^{3}}\right )+
\frac{2-r}{2(1-r)^{2}}\left (\frac{r^2(3-r)}{(1-r)^{3}}\right )+
\frac{3r-2r^{2}}{2(1-r)^{2}}\left ( \frac{r(2-r)}{(1-r)^2}\right )\\
&=& \frac{8r^{2}-10r^{3}+5r^{4}-r^{5}}{(1-r)^{5}}\\
&=& 1-  \frac{2r^{2}-5r+1}{(1-r)^{5}}
\eeqq
which implies that $\left|U_{F}(z)\right|<1$ provided $2r^{2}-5r+1>0$, i.e. $F\in {\mathcal U}$ for $|z|<r_8$,
where $r_8=\frac{5-\sqrt{17}}{4}$ is the root of the equation $2r^{2}-5r+1=0$ in the interval $(0,1)$. As for the sharpness is
concerned, we choose the univalent functions
$$f(z)=\frac{z(2-z)}{2(1-z)^{2}} ~\mbox{ and }~ g(z)=\frac{z}{(1-z)^{2}}
$$
so that $a_{n} =(n+1)/2$ and $b_{n}=n$ for $n\geq 2$. Moreover, $F$ defined by \eqref{3-10eq2} then takes the form
$$F(z)=\frac{2z(1-z)^{4}}{2-z}
$$
and thus,
$$F'(z)=\frac{4(1-z)^{3}(2z^{2}-5z+1)}{\left(2-z\right)^{2}}.
$$
It follows that $F'(r_{8})=0$ which means that $F$ is not univalent in a disk with bigger radius.

\textbf{(c)} As in the other cases,  the  assumption that $|a_{n}|\leq (n+1)/2$ and $|b_{n}|\leq 2/n$ for $n\geq 2$,  clearly gives the estimate
\beqq
\left|U_{F}(z)\right|& \leq & \left ( -1-\frac{2}{r}\log (1-r) \right )\frac{r^2(2-r)}{(1-r)^{3}}+
\frac{2-r}{2(1-r)^{2}}\left ( \frac{4-2r}{1-r} +\frac{4}{r}\log (1-r) \right ) \\
&& +\frac{3r-2r^{2}}{2(1-r)^{2}}\left ( -2-\frac{2}{r}\log (1-r) \right )\\
&=& 1-\frac{1}{(1-r)^{3}}\left ((r-3)(1-r) -\frac{4-9r+3r^2}{r}\log (1-r) \right )
\eeqq
which is less than $1$ provided $\psi (r)>0$, where
$$\psi (r)=(r-3)(1-r) -\frac{4-9r+3r^2}{r}\log (1-r),
$$
i.e. if $0<r<r_9$. This means that $F\in {\mathcal U}$ for $|z|<r_9$ where $r_9$ is as in the statement,
namely, the root of the equation $\psi (r)=0$ in the interval $(0,1)$ (see \eqref{eq10e}).

In order to prove the sharpness of the radius, we consider the univalent functions
$$f(z)=\frac{z(2-z)}{2(1-z)^{2}} ~\mbox{ and }~ g(z)=-z -2\log (1-z).
$$
We see that  $a_{n} =(n+1)/2$ and $b_{n}=2/n$ for $n\geq 2$. Using these functions, the corresponding
$F$ has the form
$$F(z)=\frac{z(1-z)^{2}}{ \left (1- (z/2)\right ) \left (-1- (2/z)\log (1-z)\right )}
$$
and thus, after some computation, one can easily see that
$$F'(z)=\frac{(1-z) [ (z-3)(1-z) - ((4-9z+3z^2)/z)\log (1-z) ]}{\left [ \left (1- (z/2)\right ) \left (-1- (2/z)\log (1-z)\right )\right]^{2}}.
$$
It follows that $F'(r_9)=0$ which means that $F$ is not univalent in a disk with bigger radius.
\hfill $\Box$


\section{Applications}\label{sec5}

There is a vast literature dedicated to the mapping properties of special functions (see for example \cite{Po97,Po98,ponu2} on Gaussian
hypergeometric functions and \cite{baricz-publ,lecture,bsk,samy-baricz} on Bessel functions).
Let us consider the normalized Bessel function $f_{\nu}:\mathbb{D}\to\mathbb{C},$ defined by
$$f_{\nu}(z)=2^{\nu}\Gamma(\nu+1)z^{1-\frac{\nu}{2}}J_{\nu}(\sqrt{z})=z+\sum_{n\geq 2}\frac{(-1)^{n-1}\Gamma(\nu+1)}{4^{n-1}(n-1)!\Gamma(\nu+n)}z^n,
$$
where $J_{\nu}$ stands for the Bessel function of the first kind. Observe that the inequality
$$\left|\frac{(-1)^{n-1}\Gamma(\nu+1)}{4^{n-1}(n-1)!\Gamma(\nu+n)}\right|\leq 1$$
is equivalent to
\begin{equation}\label{eq*}
4^{n-1}(n-1)!(\nu+1)(\nu+2)\cdots(\nu+n-1)\geq 1,
\end{equation}
and this holds for each $n\geq 2$ if $\nu\geq -\frac{3}{4}.$ Here we used that the left-hand side of the inequality \eqref{eq*} is increasing in $n$ when $\nu>-1.$ Similarly, the inequality
$$\left|\frac{(-1)^{n-1}\Gamma(\nu+1)}{4^{n-1}(n-1)!\Gamma(\nu+n)}\right|\leq n$$
is equivalent to
$$4^{n-1}n!(\nu+1)(\nu+2)\cdots(\nu+n-1)\geq 1,$$
and this holds for each $n\geq 2$ when $\nu\geq -\frac{7}{8}.$
Thus, if we consider the function $F_{\nu,\mu}:\mathbb{D}\to\mathbb{C},$ defined by
$$F_{\nu,\mu}(z)=\frac{z^3}{f_{\nu}(z)f_{\mu}(z)}=\frac{1}{2^{\nu+\mu}\Gamma(\nu+1)\Gamma(\mu+1)}\cdot \frac{z^{1+\frac{\nu+\mu}{2}}}{J_{\nu}(\sqrt{z})J_{\mu}(\sqrt{z})},$$
in view of Theorem 1 we have the following result.

\bcor\label{cor4}
The following assertions are true:
\begin{enumerate}
\item[\bf (a)] If $\nu,\mu\geq-\frac{3}{4},$ then $F_{\nu,\mu}\in\mathcal{U}$ in the disk $|z|<\frac{1}{3}.$
\item[\bf (b)] If $\nu\geq-\frac{7}{8}$ and $\mu\geq -\frac{3}{4},$ then $F_{\nu,\mu}\in\mathcal{U}$ in the disk $|z|<\frac{1}{4}.$
\item[\bf (c)] If $\nu,\mu\geq-\frac{7}{8},$ then $F_{\nu,\mu}\in\mathcal{U}$ in the disk $|z|<\frac{1}{5}.$
\end{enumerate}
\ecor

Moreover, by following the proof of \cite[Theorem 1.6]{barsza}, it can be shown that $f_{\nu}\in\mathcal{C}(-1/2)$ if and only if $\nu\geq \nu^{\star},$
where $\nu^{\star}\simeq -0.287872\dots$ is the unique root of the equation
$(2\nu-5)J_{\nu+1}(1)+5J_{\nu}(1)=0.$ Taking into account this result, Theorem \ref{3-13th2} yields the following consequence.

\bcor\label{cor5}
The following assertions are true:
\begin{enumerate}
\item[\bf (a)] If $\nu\geq\nu^{\star}$ and $\mu\geq -\frac{3}{4},$ then $F_{\nu,\mu}\in\mathcal{U}$ in the disk $|z|<r_7.$
\item[\bf (b)] If $\nu\geq\nu^{\star}$ and $\mu\geq -\frac{7}{8},$ then $F_{\nu,\mu}\in\mathcal{U}$ in the disk $|z|<r_8.$
\end{enumerate}
\ecor

Clearly, Corollaries \ref{cor4} and \ref{cor5} could be used to generate new set of functions belonging to the class $\mathcal{U}$.

\section{Concluding Remarks}\label{sec6}
Let us now consider situations dealing with the inequalities \eqref{eq4}.
A particular consideration in \eqref{eq-1} will be the five convex functions $g$, where
$$g(z)\in \left \{z,~~ z/(1 -z),~~  - \log(1 - z),~~\frac{1}{2}\log \left (\frac{1+z}{1-z}\right ),~~
-\frac{i}{\sqrt{3}}\log \left (\frac{1+\alpha z}{1+\overline{\alpha}z}\right )\right\}
$$
(where $\alpha =(-1+i\sqrt{3})/2$) so that
$$1/g'(z)\in \left \{1,~~  (1 -z)^2,~~  1-z,~~ 1-z^2,~~  1-z+z^2\right\}
$$
which led to introduce
\beqq
{\mathcal K}_1 &=&  \{f\in{\mathcal A}:\, {\rm Re\,}f'(z)>0 ~\mbox{ for $z\in \ID$}\}\\
{\mathcal K}_2 &=&  \{f\in{\mathcal A}:\, {\rm Re\,}(1-z)^2f'(z)>0 ~\mbox{ for $z\in \ID$}\}\\
{\mathcal K}_3 &=&  \{f\in{\mathcal A}:\, {\rm Re\,}(1-z)f'(z)>0 ~\mbox{ for $z\in \ID$}\}\\
{\mathcal K}_4 &=&  \{f\in{\mathcal A}:\, {\rm Re\,}(1-z^2)f'(z)>0 ~\mbox{ for $z\in \ID$}\}, ~\mbox{ and}\\
{\mathcal K}_5 &=&  \{f\in{\mathcal A}:\, {\rm Re\,}(1-z+z^2)f'(z)>0~\mbox{ for $z\in \ID$}\}.
\eeqq
Functions in ${\mathcal K}_j$ ($j=1,2,\ldots ,5$) are obviously close-to-convex in $\ID$ and each of these classes plays
a special role in certain circumstances (see \cite{MacG69, Po97,ponu2}). Moreover, these results are generalized in many ways and a necessary coefficient
condition for a function to be in ${\mathcal K}_j$ ($j=1,2,\ldots ,5)$ will be of the form \eqref{eq4}. For example, if
$${\mathcal K}_1(\alpha) = \{f\in{\mathcal A}:\, {\rm Re\,}f'(z)>\alpha~\mbox{ for $z\in \ID$}\}
$$
for some $\alpha \in [0,1)$, then $|a_n|\leq 2(1-\alpha)/n$ for all $n\geq 2$.  It is important to point out that functions in
${\mathcal K}_1:={\mathcal K}_1(0)\equiv \mathcal{R}$ are not necessarily starlike.

As remarked in the introduction and also from the statement of our main results, one can observe from a careful analysis that the results
and the methodology of this article are
applicable in general settings although we deal in this paper only important cases. In order to have this feeling,
it would be interesting to recall also the following few cases:
\bee
\item[(a)] if $f\prec g$ and $g\in \mathcal{R}\cup {\mathcal C}$, then $|a_n|\leq 1$ for all $n\geq 2$ (See \cite[Theorem 5]{Hall74}
and \cite[p.~195, Theorem 6.4]{Du})
 \item[(b)] if $f\prec g$ and $g\in {\mathcal S}^*$, then $|a_n|\leq n$ for all $n\geq 2$ (See \cite[p.~195, Theorem 6.4]{Du})
\item[(c)] if $f\prec g$ and $g\in {\mathcal S}$, then $|a_n|\leq n$ for all $n\geq 2$ (See \cite[p.~196]{Du} and \cite{DeB1})
\item[(d)] if $f\in {\mathcal A}$ with ${\rm Re}\, (f(z)/z)>\alpha$ in $\ID$
for some $\alpha \in [0,1)$, then $|a_n|\leq 2(1-\alpha)$ for all $n\geq 2$
\item[(e)] if $f\in {\mathcal G}$, then $|a_n|\leq 1/(n(n-1))$ for all $n\geq 2$, where
$$ {\mathcal G} = \left \{f\in {\mathcal A}:\,{\rm Re\,}\left(1+\frac{zf''(z)}{f'(z)}\right)<\frac{3}{2}~\mbox{ for $z\in \ID$} \right \}.
$$
Each $f\in{\mathcal G}$ is known to be starlike in $\ID$. See \cite{ObPoWi2013} for details and further investigation on this class.
\eee

Now, it is worth pointing out that our results are applicable for the cases (a)--(c). The cases (d)--(e) along with other general assumptions could
be used to  prove several new results. In addition to this remark, it is worth recalling some related problems in the literature for which our methodology works,
because there are many other important classes of functions for which our results could be stated after some investigations.
A sequence of real numbers $\{c_k\}_{k\geq 1}$
is said to be totally monotonic (resp. $n$-times monotonic) if
$\Delta ^0 c_k =c_k\geq 0$ and $\Delta ^m c_k=\Delta ^{m-1} c_k-\Delta ^{m-1} c_{k+1}\geq 0$ for all $m\geq 1$, $k\geq 1$
(resp.  $\Delta ^m c_k\geq 0$ for $m=0,1, \ldots, n$). An $n$-times monotonic sequence is denoted by ${\mathcal M}_n$.

In \cite{wirths2}, Wirths has shown that if $f\in {\mathcal A}$ and $\{a_k\}_{k\geq 1}$ are totally monotonic, then the largest
number $r^*$ such that $f\in {\mathcal S}^*$ in $ |z|<r^*$ is   $r^*\approx 0.934$.  Similar results are obtained for other related families
of functions from ${\mathcal A}$. In another article \cite{wirths3}, Wirths has shown that if $ \{(k+1-2\beta)a_k\}_{k\geq 1}\in {\mathcal M}_2$ for
some $\beta \in [0,1)$, then $f\in {\mathcal S}^*(\beta)$, where ${\mathcal S}^*(\beta)$ denotes the class of starlike functions of order
$\beta$ in $\ID$. Finally, a necessary and sufficient condition for the coefficients of $f\in {\mathcal A}$ to be totally monotonic is that
$$f(z)=\int_0^1\frac{z}{1-tz}\,d\mu(t)
$$
for some probability measure $\mu(t)$ defined on the unit interval $[0, 1]$. Such functions are known to
be univalent in $\ID$. Thus, our results are applicable for analytic functions in $\ID$ whose Taylor coefficients form a totally monotonic sequence.
Moreover, Wirths \cite{wirths1} has found the radius of starlikeness for the class of functions of the above form,
and the radius of convexity for  this class was proved to be
$1/\sqrt{2}$ (see \cite{SilSilTel-78}).

\subsection*{Acknowledgements}
The work of the second author was supported by MNZZS Grant, No. ON174017, Serbia.
The third author is on leave from the Indian Institute of Technology Madras, India.

\end{document}